\documentclass[11pt]{amsart}%
\usepackage{amsmath,amssymb,latexsym,graphicx}
\usepackage{enumerate}
\usepackage{amsmath}
\usepackage{amsfonts}
\usepackage{amssymb}
\usepackage{color}
\usepackage{cite}
\usepackage{graphicx}%
\providecommand{\U}[1]{\protect\rule{.1in}{.1in}}
\newtheorem{theorem}{Theorem}[section]
\newtheorem{lemma}[theorem]{Lemma}

\definecolor{gree}{rgb}{0.,0.4,0}
\begin{document}
\title[On the exactness of the universal backprojection formula]{On the exactness of the universal backprojection formula for the spherical means
Radon transform}
\author{Mark~ Agranovsky and Leonid~ Kunyansky}

\begin{abstract}
The spherical means Radon transform $\mathcal{M}f(x,r)$ is defined by the
integral of a function $f$ in $\mathbb{R}^{n}$ over the sphere $S(x,r)$ of
radius $r$ centered at a $x$, normalized by the area of the sphere. The
problem of reconstructing $f$ from the data $\mathcal{M}f(x,r)$ where $x$
belongs to a hypersurface $\Gamma\subset\mathbb{R}^{n}$ and $r \in(0,\infty)$
has important applications in modern imaging modalities, such as photo- and
thermo- acoustic tomography. When $\Gamma$ coincides with the boundary
$\partial\Omega$ of a bounded (convex) domain $\Omega\subset\mathbb{R}^{n}$, a
function supported within $\Omega$ can be uniquely recovered from its
spherical means known on $\Gamma$. We are interested in explicit inversion
formulas for such a reconstruction.

If $\Gamma=\partial\Omega$, such formulas are only known for the case when
$\Gamma$ is an ellipsoid (or one of its partial cases). This gives rise to the
natural question: can explicit inversion formulas be found for other closed
hypersurfaces $\Gamma$? In this article we prove, for the so-called "universal
backprojection inversion formulas", that their extension to non-ellipsoidal
domains $\Omega$ is impossible, and therefore ellipsoids constitute the
largest class of closed convex hypersurfaces for which such formulas hold.

\end{abstract}
\maketitle

\emph{Keywords:} Universal backprojection formula, thermoacoustic tomography, explicit inversion formula, spherical means

\section{Introduction}

\subsection{Formulation of the problem and the main result}

Given a continuous function $f$ in $\mathbb{R}^{n}$ define by $\mathcal{M}%
f(x,r),\ x\in\Gamma,\ r>0$ the spherical mean
\[
\mathcal{M}f(x,r)=\frac{1}{|\mathbb{S}^{n-1}|}\int\limits_{\mathbb{S}^{n-1}%
}f(x+r\theta)\,dS(\theta),
\]
where $\mathbb{S}^{n-1}$ is the unit sphere in $\mathbb{R}^{n},$
$|\mathbb{S}^{n-1}|$ is the $(n-1)$-dimensional area of $\mathbb{S}^{n-1},$
$dS$ is the surface Lebesgue measure on $\mathbb{S}^{n-1}.$

The spherical mean operator takes functions $f$ in $\mathbb{R}^{n}$ to
functions $\mathcal{M}f(x,r)$ on $\mathbb{R}^{n}\times\lbrack0,\infty]$. The
problem of reconstructing $f$ from the spherical means $\mathcal{\mathcal{M}%
}f(x,r)$, with centers $x$ located on a certain hypersurface $\Gamma$ arises
in inverse problems and modern modalities of imaging, such as thermo- and
photoacoustuc tomography (TAT/PAT)
\cite{Kruger95,Oraev94,Kruger99}. There is an extensive literature on the
analytical, computational and applied aspects of this subject; we refer the
reader to surveys \cite{KuKu,KuKuHand} and references therein.

The important case we consider here is when the centers $x$ belong to a closed
hypersurface $\Gamma$ which coincides with the boundary $\partial\Omega$ of a
bounded domain $\Omega$ that, in turn, contains the support of function $f$.
The transform $\mathcal{M}$ is invertible in this case
\cite{Finch04,Finch07,US}, but explicit inversion formulas are only known for
ellipsoidal surfaces $\Gamma$ (including spheres, spheroids, etc.)
\cite{Finch04,Finch07,FR3,Haltm-ellipsoid,Haltm-universal-SIAM,Kun-explicit,Nat12,Pala,Pal-book-2016,Salman}%
.

In this paper we focus on a certain type of inversion formulas, first proposed
(in an equivalent form) in \cite{Wang-sphere} for the partial case when
$\partial E$ is a sphere in $\mathbb{R}^{3}$ and later extended to spheres in
arbitrary dimensions and, finally, to arbitrary ellipsoids in $\mathbb{R}^{n}$
\cite{Haltm-ellipsoid,Haltm-universal-JMAA,Haltm-universal-SIAM,Nat12}. In the
latter case, when $\Omega$ is an ellipsoidal domain $E$ in $\mathbb{R}^{n}$
bounded by an ellipsoid $\Gamma=\partial E$ the inversion formula for
reconstruction function $f$ from the values of its spherical means
$\mathcal{M} f(x,r), \ x \in\partial E, r \geq0,$  has the
following form: \cite{Haltm-universal-JMAA}:
\begin{equation}
f(x_{0})=\left[  \mathcal{N}_{E}\mathcal{M}f\right]  (x_{0})=\frac
{(-1)^{\frac{n-2}{2}}|\mathbb{S}^{n-1}|}{2\pi^{n}}\int\limits_{\partial E}%
\nu_{x}\cdot(x_{0}-x)\int\limits_{0}^{\infty}\frac{\left[  \partial
_{r}\mathcal{D}_{r}^{n-2}r^{n-2}(\mathcal{M}f)\right]  (x,r)}{r^{2}%
-|x_{0}-x|^{2}|}\,dr\,dS(x), \label{E:inv-even}%
\end{equation}
(when $n\geq2$ is even), and
\begin{equation}
f(x_{0})=\left[  \mathcal{N}_{E}\mathcal{M}f\right]  (x_{0})=\frac
{(-1)^{\frac{n-3}{2}}|\mathbb{S}^{n-1}|}{4\pi^{n-1}}\int\limits_{\partial
E}\nu_{x}\cdot\frac{x_{0}-x}{|x_{0}-x|}\left[  \partial_{r}\mathcal{D}%
_{r}^{n-2}r^{n-2}(\mathcal{M}f)\right]  (x,|x_{0}-x|))\,dS(x),
\label{E:inv-odd}%
\end{equation}
(when $n\geq3$ is odd). Here $f\in C^{\infty}$ with $\mathrm{supp}\ f\subset
E,$ \ $x_{0}$ is an arbitrary point in $E,$ $\nu_{x}$ denotes the exterior
unit normal to $\partial E$, and $\mathcal{D}_{r}=(2r)^{-1}\partial_{r}$ is
the operator of differentiation with respect to $r^{2}.$ By $\mathcal{N}_{E},$
we denote the operator, of backprojection type, in the right hand side,
applied to $\mathcal{M}f.$ In operator terms, it is the left inverse operator
to the restricted spherical means Radon transform $\mathcal{M} f(x,r)\vert
_{\Gamma\times[0,\infty)}.$ Versions of the inversion formulas for
 certain unbounded quadratic surfaces were also obtained
\cite{Wang-universal,Wang-cylinder,Wang-various}.

Thus, among bounded closed convex observation surfaces $\Gamma,$ ellipsoids
are the only those for which explicit inversion formulas are known by now. The
natural question arises, whether nice inversion formulas can be constructed
for other surfaces? Specifically, we are interested in answer to this question
addressed inversion formulas (\ref{E:inv-even}), (\ref{E:inv-odd}):
\textit{are these formulas true for bounded convex smooth hypersurfaces
}$\Gamma$\textit{ other than ellipsoids }$\partial E$\textit{? }The goal of
this paper is to show that the answer is \textit{negative:} ellipsoids
constitute the largest class of smooth bounded convex hypersurfaces $\Gamma$
for which these inversion formulas are valid. More precisely, the following
theorem holds:

\begin{theorem}
\label{T:Main} Let $\Omega$ be a bounded convex domain in $\mathbb{R}%
^{n},\ n\geq2,$ and $\Gamma=\partial\Omega\in C^{\infty}.$ Suppose that the
inversion formula (\ref{E:inv-even}) or (\ref{E:inv-odd}), where $\partial E$
is replaced by $\partial\Omega,$ holds for all $x_{0}\in\Omega$ and for any
function $f\in C^{\infty}(\mathbb{R}^{n})$ supported in $\Omega.$ Then the
hypersurface $\partial\Omega$ is an ellipsoid $\partial E.$
\end{theorem}

\subsection{"Universal backprojection formulas"}

Inversion formulas (\ref{E:inv-even}), (\ref{E:inv-odd}) and
their equivalents have received the name of "universal backprojection
formulas" \cite{Wang-universal,Haltm-universal-SIAM,Haltm-universal-JMAA}. Our
proof of Theorem \ref{T:Main} relies on the most general form of these
formulas given in \cite{Haltm-universal-SIAM}. Namely, it has been found that
if in the right hand side of formulas (\ref{E:inv-even}), (\ref{E:inv-odd})
the ellipsoid $\partial E$ is replaced by the boundary $\Gamma=\partial\Omega$
of an arbitrary bounded convex smooth domain, then an additional "error" term
$\mathcal{K}_{\Omega}f$ will appear:
\begin{equation}
f(x_{0})=\left[  \mathcal{N}_{\Omega}\mathcal{M}f\right]  (x_{0})+\left[
\mathcal{K}_{\Omega}f\right]  (x_{0}), \label{E:univ}%
\end{equation}
where $\mathcal{N}_{\Omega}$ is given by (\ref{E:inv-even}), (\ref{E:inv-odd})
( with $\Omega$ in place of $E$) and $\mathcal{K}_{\Omega}$ is the integral
operator
\begin{equation}
\left[  \mathcal{K}_{\Omega}f\right]  (x_{0})=\int\limits_{\Omega}k_{\Omega
}(x_{0},x_{1})f(x_{1}) \, dx_{1}, \label{E:K-kernel}%
\end{equation}
where the kernel $k_{\Omega}$ associated with $\Omega$ is explicitly given by
formulas (\ref{E:k-even}), (\ref{E:k-odd}) in Section \ref{S:error}.

Of course, the name "universal backprojection formula" for (\ref{E:univ}) is
somewhat  misleading, since  extending formula
(\ref{E:univ})
to a larger class of the domains $\Omega$, in general, results in the loss of
its inversion property. Indeed, equation (\ref{E:univ}) leads to an integral
equation of the second type $(I-\mathcal{K}_{\Omega})f=N_{\Omega}\left(
\mathcal{M}f\right)  $ for the unknown function $f$ rather than to an explicit
expression for $f.$ Therefore, in order to have a true inversion formula
$f=\mathcal{N}_{\Omega}\left(  \mathcal{M}f\right)  $  one needs
to guarantee that the error term vanishes. Thus, we would like to
characterize all domains $\Omega$ for which $\mathcal{K}_{\Omega}=0.$ In the latter case
we will say that the universal inversion formula (\ref{E:univ}) is
\emph{exact}. In these terms, Theorem \ref{T:Main} can be translated as
follows: non-ellipsoidal domains necessarily produce the non-zero "error" term
$\mathcal{K}_{\Omega}$ and hence "universal inversion formula" (\ref{E:univ})
is exact \textbf{if and only if} the boundary of $\Omega$ is an ellipsoid.

The rest of this paper is arranged as follows. In the next section we define
the integral transforms needed to properly present the explicit expression for
$\mathcal{K}_{\Omega}$ and prove two lemmas necessary for the further
exposition. Theorem \ref{T:Main} is proven in Section \ref{S:Main}. We
conclude with the further discussion of our results in Section
\ref{S:discussion}.

\section{Preliminaries\label{S:Lemmas}}

\subsection{Radon and Hilbert transforms}

Below we recall several well known facts about the Radon and Hilbert transforms.

The Radon transform $\mathcal{R}f$ of a compactly supported smooth function
$f$ is defined\cite{natt-book-1} as
\[
\lbrack\mathcal{R}f](\omega,p)= \int\limits_{\Pi(\omega,p)}f(x)\,dA(x),\qquad
(\omega,p)\in\mathbb{S}^{n-1}\times\mathbb{R},
\]
where $\mathbb{S}^{n-1}$ is the unit sphere in $\mathbb{R}^{n}$, $\Pi
(\omega,p)$ is the hyperplane defined by the equation $x\cdot\omega=p$, and
$dA(x)$ is the standard measure on $\Pi(\omega,p)$. Obviously, $[\mathcal{R}%
f](\omega,p)$ vanishes for all $(\omega,p)$ such that $\Pi(\omega,p)$ does not
intersect the support $\Omega$ of $f$.

The Hilbert transform $\mathcal{H}F$ of a smooth, sufficiently fast decaying
function $F(t),$ $t\in\mathbb{R}$, is defined by the following formula
\cite{King}%
\[
\lbrack\mathcal{H}F](t) = \frac{1}{\pi}p.v.\int\limits_{\mathbb{R}}\frac
{F(s)}{t-s}\,ds\text{,}%
\]
where $p.v.$ stands for the principal value of the integral. The Hilbert
transform can be extended to less smooth functions and distributions by
continuity. It is self-invertible; more precisely $\mathcal{H(H}F)=-F.$ The
following intertwining relation holds (Section 4.7,\cite{King}):
\begin{equation}
\left[  \mathcal{H}(s\varphi(s))\right]  (t)=t\left[  \mathcal{H}%
\varphi\right]  (t)-\frac{1}{\pi}\int\limits_{\mathbb{R}}\varphi(t) \, dt\text{.}
\label{E:intertwine}%
\end{equation}
Let $\chi_{\lbrack a,b]}(s)$ be the characteristic function of the interval
$[a,b].$ The Hilbert transform of the function $\chi_{\lbrack-1,1]}%
(s)\sqrt{(1-s)(1+s)}$ is well-known (formula 11.343,\cite{King}):%
\[
\left[  \mathcal{H}\left(  \chi_{\lbrack-1,1]}(s)\sqrt{(1-s)(1+s)}\right)
\right]  (t)=t,\quad t\in\lbrack a,b].
\]
By a linear change of variables one obtains the Hilbert transform of
$\chi_{\lbrack a,b]}(s)\sqrt{(b-s)(s-a)}$:%

\begin{equation}
\left[  \mathcal{H}\left(  \chi_{\lbrack a,b]}(s)\sqrt{(b-s)(s-a)}\right)
\right]  (t)=t-\frac{b+a}{2},\quad t\in\lbrack a,b]. \label{E:Hilbert_arch}%
\end{equation}
We will also make use of the existence of the so-called finite inverse Hilbert
transform \cite{you2006explicit}. Namely, if a continuous function $F(t)$ is
supported on an interval $[a,b]$, then
\begin{equation}
F(t)=\frac{1}{\pi\sqrt{(b-t)(t-a)}}\left(  \int\limits_{a}^{b}\frac
{[\mathcal{H}F](s)}{s-t}\sqrt{(b-s)(s-a)} \, ds+\int\limits_{a}^{b}F(s) \, ds\right)
,\ t\in\lbrack a,b]. \label{E:finite_Hilbert}%
\end{equation}

\subsection{The error operator}

\label{S:error}

Our analysis of the universal backprojection formula is based on the results
of \cite{Nat12,Haltm-convex,Haltm-universal-SIAM}, that give explicit
expressions for the error arising when this formula is used under the
following assumptions: $f(x)\ $is is compactly supported strictly inside of a
bounded (strictly) convex open domain $\Omega\subset\mathbb{R}^{n}$ with an
infinitely smooth boundary $\partial\Omega$ and the measuring surface $\Gamma$
coincides with $\partial\Omega.$ Then the universal backprojection operator
$\mathcal{N}_{\Omega}$, when applied to spherical means $\mathcal{M}f$ will
produce the error $\mathcal{K}_{\Omega}f$, see equation (\ref{E:univ}), with
operator $\mathcal{K}_{\Omega}$ in the form (\ref{E:K-kernel}). The latter
kernel, according to \cite{Nat12,Haltm-convex,Haltm-universal-SIAM}, has the
following form:
\begin{equation}
k_{\Omega}(x_{0},x_{1})=c\ \frac{\partial_{s}^{n}\mathcal{[H(R}\chi_{\Omega
})](\omega_{\ast}(x_{0},x_{1}),s_{\ast}(x_{0},x_{1}))}{|x_{0}-x_{1}|^{n-1}},
\label{E:k-even}%
\end{equation}
if $n$ is even, and
\begin{equation}
k_{\Omega}(x_{0},x_{1})=c\ \frac{\partial_{s}^{n}[\mathcal{R}\chi_{\Omega
}](\omega_{\ast}(x_{0},x_{1}),s_{\ast}(x_{0},x_{1}))}{|x_{0}-x_{1}|^{n-1}},
\label{E:k-odd}%
\end{equation}
if $n$ is odd. Here $\mathcal{R}\chi_{\Omega}$ is the Radon transform of the
characteristic function $\chi_{\Omega}$ of $\Omega,$ the Hilbert transform
$\mathcal{H}$ acts with respect to the second variable of the pair
$(\omega,s)$, and functions $\omega_{\ast}(x_{0},x_{1})$ and $s_{\ast}%
(x_{0},x_{1})$ are defined as follows:%
\[
\omega_{\ast}(x_{0},x_{1})=\frac{x_{0}-x_{1}}{|x_{0}-x_{1}|},\ s_{\ast}%
(x_{0},x_{1})=\frac{1}{2}\frac{|x_{1}|^{2}-|x_{0}|^{2}}{|x_{1}-x_{0}|}.
\]
It has been proven that if the domain $\Omega$ is an ellipsoid then the error
operator $\mathcal{K}_{\Omega}$ vanishes
\cite{Nat12,Haltm-convex,Haltm-universal-SIAM}, and hence the expression
$\mathcal{N}_{\Omega}(\mathcal{M}f)$ represents the exact inversion and
returns $f$. As it was already mentioned, the main result of the present paper
( Theorem \ref{T:Main}) is that the converse statement is true: the universal
backprojection inversion formula (\ref{E:univ})
is exact only for ellipsoids.


We proceed with two lemmas which we will need in the proof of Theorem
\ref{T:Main}.

\subsection{Two lemmas}

\begin{lemma}
\label{L:Hilbert} Let $[a,b]$ be a segment on the real line and let $F\in
C(\mathbb{R})$ be supported in the segment $[a,b]$ . If there exists a
polynomial $P(t)$ such that $\left[  \mathcal{H}F\right]  (t)=P(t)$ for all
$t\in\lbrack a,b]$ then
\begin{equation}
F(t)=\frac{Q(t)}{\sqrt{(b-t)(t-a)}},\ t\in\lbrack a,b], \label{E:PQ}%
\end{equation}
where $Q(t)$ is a polynomial of the degree $\deg Q\leq\deg P+1.$
\end{lemma}

\begin{proof}
Consider the Hilbert transform of a function $\chi_{\lbrack a,b]}(s)s^{k}%
\sqrt{(b-s)(s-a)}$ with integer $k.$ If $k=0$, identity (\ref{E:Hilbert_arch})
yields%
\[
\mathcal{H}\left(  \chi_{\lbrack a,b]}(s)\sqrt{(b-s)(a-s)}\right)
(t)=t+c_{0},\quad t\in\lbrack a,b]\text{.}%
\]
with some constant $c_{0}.$ For $k>0$ equation (\ref{E:intertwine}) leads to%
\[
\mathcal{H}\left(  \chi_{\lbrack a,b]}(s)s^{k}\sqrt{(b-s)(a-s)}\right)
(t)=t\left\{  \mathcal{H}\left(  \chi_{\lbrack a,b]}(s)s^{k-1}\sqrt
{(b-s)(a-s)}\right)  (t)+c_{k}\right\}  ,
\]
where $c_{k}$ is yet another constant. By induction, the above two equations
imply that \newline$\mathcal{H}\left(  \chi_{\lbrack a,b]}(s)s^{k}%
\sqrt{(b-s)(a-s)}\right)  $ is a polynomial of degree $k+1$. Thus, if $P(s)$
is a polynomial of degree $\deg P,$ then $\mathcal{H}\left(  \chi_{\lbrack
a,b]}(s)P(s)\sqrt{(b-s)(s-a)}\right)  (t)$ is a polynomial of degree $\deg
P+1.$ On the other hand, $(\mathcal{H}F)(t)=P(t)$ for $t\in\lbrack a,b]$ and
hence formula (\ref{E:finite_Hilbert}) reads as%
\[
F(t)=\frac{1}{\sqrt{(b-t)(t-a)}}\left(  -\left[  \mathcal{H}\left(
\chi_{\lbrack a,b]}(s)P(s)\sqrt{(b-s)(s-a)}\right)  \right]  (t)+\int_{a}%
^{b}F(s) \, ds\right)  ,\ t\in\lbrack a,b],
\]
thus proving the lemma.
\end{proof}

Given a unit vector $\omega\in\mathbb{S}^{n-1},$ define%

\begin{align}
&  h_{\Omega}^{+}(\omega)=\sup_{ x \in\Omega} x \cdot\omega,\\
&  h_{\Omega}^{-}(\omega)=\inf_{ x \in\Omega} x \cdot\omega.
\end{align}

\begin{figure}[t]
\begin{center}
\includegraphics[scale = 0.8]{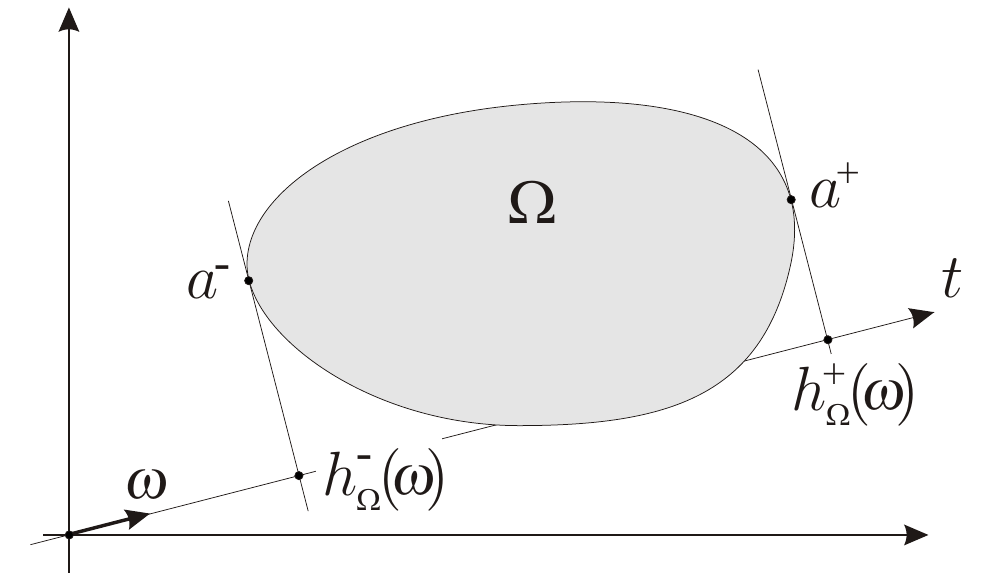}
\end{center}
\caption{Geometric meaning of functions $h_{\Omega}^{+}(\omega)$ and
$h_{\Omega}^{-}(\omega)$ }%
\label{F:geometry}%
\end{figure}

The function $h_{\Omega}(\omega)=h_{\Omega}^{+}(\omega),$ $\omega\in
\mathbb{S}^{n-1}$, is called \textbf{the support function} of the domain
$\Omega.$

The functions $h_{\Omega}^{\pm}$ are related by the formula
\[
h_{\Omega}^{-}(\omega)=-h_{\Omega}(-\omega),\qquad\omega\in\mathbb{S}^{n-1}.
\]
In the case of ellipsoidal domains, the support function is the square root of
a quadratic polynomial. For example, for the domain $E$ bounded by the
ellipsoid
\[
\partial E=\left\{  \sum_{j=1}^{n}\frac{x_{j}^{2}}{a_{j}^{2}} = 1\right\}
\]
we have
\[
h_{E}(\omega)=\sqrt{\sum_{j=1}^{n}a_{j}^{2}\omega_{j}^{2}}.
\]

A hyperplane $\{x\cdot\omega=t\}$ meets the domain $\Omega$ if and only if
$h_{\Omega}^{-}(\omega)<t<h_{\Omega}(\omega).$ The limit cases $t=h_{\Omega
}^{\pm}(\omega)$ correspond to the tangent hyperplanes to $\partial\Omega$ at
the points $a^{\pm}\in\partial\Omega$ where the exterior unit normal vectors
are $\nu_{\partial\Omega}(a^{\pm})=\pm\omega$, as illustrated in Figure
\ref{F:geometry}.

The behavior of the Radon transform $\left[  \mathcal{R}\chi_{\Omega}\right]
(\omega,t)$ near the tangent planes is given by the following Lemma.

\begin{lemma}
\label{L:zeros} For a dense set of the direction vectors $\omega\in
\mathbb{S}^{n-1},$ the following asymptotic relation holds with some nonzero
constants $c^{\pm}$:
\begin{equation}
\big[\mathcal{R}\chi_{\Omega}\big](\omega,t)=c^{\pm}(t-h_{\Omega}^{\pm}%
(\omega))^{\frac{n-1}{2}}(1+\mathit{o}(1)),\ t\rightarrow h_{\Omega}^{\pm
}(\omega)\mp0. \label{E:asymp}%
\end{equation}

\end{lemma}

\begin{proof}
We will use the notation $\Gamma=\partial\Omega.$ The hypersurface $\Gamma$ is
infinitely differentiable. Let $\kappa_{\Gamma} (a), \ a \in\Gamma$ be the
Gaussian curvature, i.e. the product of the principal curvatures of the
$C^{\infty}$ hypersurface $\Gamma$ at the point $a.$

Denote by $\gamma$ the Gauss mapping
\[
\gamma:\Gamma\ni a\rightarrow\nu_{\Gamma}(a)\in\mathbb{S}^{n-1},
\]
which maps a point $a\in\Gamma$ to the exterior unit normal vector
$\gamma(a)=\nu_{\Gamma}(a)$ to $\Gamma$ at the point $a.$ Since $\Gamma$ is
strictly convex, $\gamma$ is a one-to-one mapping. It is differentiable and
Gaussian curvature $\kappa_{\Gamma}(a)$ equals to Jacobian determinant
$\kappa(a)=J_{\gamma}(a)$ of $\gamma$ at the point $a.$ Therefore, the points
$a$ with $\kappa_{\gamma}(a)\neq0$ (\textit{non-degenerate points}) constitute
the set $\mathrm{Reg}_{\gamma}$ of regular points of the mapping $\gamma,$
while the set of points $a$ of zero Gaussian curvature coincides with the
critical set $\mathrm{Crit}_{\gamma}.$ By Sard's theorem (see e.g., \cite{Milnor},
Section 2, p.10; Section 3, p.16) , $\gamma(\mathrm{Crit}_{\gamma})$ has the Lebesgue
measure zero on $\mathbb{S}^{n-1}$ , while the set $\gamma(\mathrm{Reg}_{\gamma})$ of
regular values is a dense subset of $\mathbb{S}^{n-1}.$ This subset consists
of the \textit{regular directions} $\omega$ which are normal vectors
$\omega=\nu_{\Gamma}(a)$ at non-degenerate points $a,$ i.e. points of nonzero
Gaussian curvature.

Let $\omega\in\mathbb{S}^{n-1}$ be a regular direction, $\omega=\nu_{\Gamma
}(a)$, $a\in\Gamma.$ Applying a suitable translation and orthogonal
transformation, we can assume that $a=0,\ \omega=(0,...,0,1).$ Then the
tangent plane $T_{a}(\Gamma)$ is the coordinate plane $x_{n}=0$ and the domain
$\Omega$ is contained in the half-space $x_{n}<0.$ In this case $h_{\Omega
}^{+}(\omega)=0.$ Moreover, after performing a suitable non-degenerate linear
transformation we can make the equation of $\Gamma,$ near $a=0,$ to be:%

\begin{equation}
\label{E:x_n=}x_{n}=-\frac{1}{2}\left(  c_{1}x_{1}^{2}+\cdots+c_{n-1}%
x_{n-1}^{2}\right)  +\mathit{o}\left(  |x^{\prime}|^{2}\right)  ,\ (x_{1}%
,...,x_{n-1})=x^{\prime}\rightarrow0.
\end{equation}
The new axes $x_{j}$, {$j=1,...,n-1$} are directions of the vectors of
principal curvatures and the coefficients $c_{j}$ are the values of the
principal curvatures at the point $a=0\in\Gamma.$ The Gaussian curvature at
$a=0$ is $\kappa_{\Gamma}(0)=c_{1}\cdots c_{n-1}$. All the applied
transformations preserve regular points, hence $\kappa_{\Gamma}(0)\neq0$.
Therefore, none of $c_{j}$ is zero and since $c_{j}\geq,0$ due to the
convexity of $\Gamma,$ we have $c_{j}>0$ for all $j.$

Since, after the above transformations, we have $\omega=(0, \cdots, 0,1)$ the
hyperplane $x\cdot\omega=t$ now is given by the equation $x_{n}=t,$ with
$t<0.$ The main term of $\mathrm{Vol}_{n-1}(\Omega\cap\{x_{n}=t\})$ near
$t=0$, is determined by the main term of expansion (\ref{E:x_n=}), i.e., by
the volume of the ellipsoid $-t=c_{1}x_{1}^{2}+\cdots+c_{n-1}x_{n-1}^{2},$
which is equal to
\[
c(-t)^{\frac{n-1}{2}};\ c=\frac{(2\pi)^{\frac{n-1}{2}}}{\Gamma(\frac{n+1}%
{2})\sqrt{\kappa_{\Gamma}(a)}}.
\]

Thus, for the specific choice $a=0$ and $\omega=(0,...0,1),$ we have the
following asymptotic formula (see e.g. \cite{GGV} Ch.1, Section 1, 7):
\[
\left[  \mathcal{R}\chi_{\Omega}\right]  (\omega,t)=\mathrm{Vol}_{n-1}%
\{x_{n}=t\}=c(-t)^{\frac{n-1}{2}}+o(|t|^{\frac{n-1}{2}}),t\rightarrow-0,
\]
Performing the inverse affine transformation, we obtain the asymptotic formula
(\ref{E:asymp}) near the point $h_{\Omega}^{+}(\omega)$ with some new nonzero
constant $c^{+}$. There are two points $a^{\pm}\in\Gamma$ with parallel
tangent planes and opposite exterior unit normal vectors $\nu_{\Gamma}(a^{\pm
})=\pm\omega,$ hence, by repeating the argument for the point on $\Gamma$ with
the exterior unit normal vector $-\omega,$ we obtain the similar asymptotic
near the point $h_{\Omega}^{-}(\omega)=-h_{\Omega}^{+}(-\omega).$ Lemma is proved.
\end{proof}

\section{Proof of Theorem \ref{T:Main}\label{S:Main}}

Let $\Omega$ be a domain in $\mathbb{R}^{n}$ satisfying all the conditions of
Theorem \ref{T:Main}. The exactness of the universal backprojection formula
for any function $f(x)$ implies that the kernel $k_{\Omega}(\omega_{\ast
}(x_{0},x_{1}),s_{\ast}(x_{0},x_{1}))$ vanishes for all $x_{0},x_{1}\in
\Omega.$ The geometric meaning of the variables $\omega_{\ast}(x_{0},x_{1})$
and $s_{\ast}(x_{0},x_{1})$ is that the hyperplane $x\cdot\omega_{\ast}%
(x_{0},x_{1})=s_{\ast}(x_{0},x_{1})$ is orthogonal to the segment
$[x_{0},x_{1}]$ and passes through its midpoint. Obviously, every hyperplane
intersecting the interior of $\Omega$ can be obtained by choosing certain
$x_{0},x_{1}\in\Omega$. Therefore, $k_{\Omega}(\omega,t)=0$ for all
$(\omega,t)\in\mathbb{S}^{n-1}\times(h_{\Omega}^{-}(\omega),h_{\Omega}%
^{+}(\omega))$. This is also trivially true for $t$ lying outside of the
interval $(h_{\Omega}^{-}(\omega),h_{\Omega}^{+}(\omega))$, so $k_{\Omega
}(\omega,t)=0$ for all $(\omega,t)\in\mathbb{S}^{n-1}\times\mathbb{R}$. This
in turn, implies that, for a fixed $\omega,$ functions $[\mathcal{H(R}%
\chi_{\Omega})](\omega,t)$ in the even dimensional case and $[\mathcal{R}%
\chi_{\Omega}](\omega,t)$ in the odd dimensional case are polynomials in $t$
of degree not exceeding $n-1$.

Below we consider the cases of even and odd $n$ separately. We start with odd dimensions.

\subsection{The case of odd $n$}

The condition $k_{\Omega}(\omega,s)=0, \ h_{\Omega}^{-}(\omega)<t< h_{\Omega
}^{+}(\omega)$ and expression (\ref{E:k-odd}) for $k_{\Omega}$ imply that
\begin{equation}
\lbrack\mathcal{R}\chi_{\Omega}](\omega,t)=P_{\omega}(t),\ h_{\Omega}^{-}
(\omega)<t< h_{\Omega}^{+}(\omega), \label{E:keven}%
\end{equation}
where
\[
P_{\omega}(t)=\sum_{k=0}^{n-1}p_{k}(\omega)t^{k}%
\]
is a polynomial of degree at most $n-1,$ with coefficients $p_{k}(\omega)$ -
continuous functions on the unit sphere $\mathbb{S}^{n-1}.$

Domains $\Omega$ with polynomial dependence of $[\mathcal{R}\chi_{\Omega
}](\omega,t)$ on $t$ are called \textit{polynomially integrable}
\cite{agranovsky2018poly}. There is no such domains in even dimensions
\cite{agranovsky2018poly}. Koldobsky, Merkurjev and Yaskin proved in
\cite{koldobsky2017} that, in odd dimensions, all polynomially integrable
domains with infinitely smooth boundaries are ellipsoids. Therefore, we could
refer here to that result. However, in our case we have an additional
information about the degree of polynomial $P_{\omega}$ which allows us to use
simpler arguments than those in \cite{koldobsky2017}. These arguments are
given in \cite{agranovsky2018poly} (see also \cite{agranovsky2019alg}). They
are based on Lemma \ref{L:zeros} and on the range conditions for the Radon
transform. We will present them here, to make the presentation self-contained
and because these arguments extend to the case of even $n,$ where the result
of \cite{koldobsky2017} is not directly applicable.

Lemma \ref{L:zeros} asserts that for a dense set of $\omega\in\mathbb{S}%
^{n-1}$ polynomial $P_{\omega}(t)$ has zeros at the points $t_{1}=h_{\Omega
}^{-}(\omega)$ and $t_{2}=h_{\Omega}^{+}(\omega),$ each zero of multiplicity
$\frac{n-1}{2}.$ This means $\partial_{t}^{\frac{n-1}{2}}P_{\omega}(t)=0$ when
$t=t_{1}$ or $t=t_{2}$ and by continuity with respect to $\omega$ these
properties extend to all $\omega\in\mathbb{S}^{n-1}.$ On the other hand, we
have the upper bound $\deg P_{\omega}\leq n-1.$ Therefore, we conclude that
polynomial $P_{\omega}$ can be represented in the form%
\begin{equation}
P_{\omega}(t)=A(\omega)(t-t_{1})^{\frac{n-1}{2}}(t_{2}-t)^{\frac{n-1}{2}%
}=A(\omega)(t-h_{\Omega}^{-}(\omega))^{\frac{n-1}{2}}(h_{\Omega}^{+}%
(\omega)-t)^{\frac{n-1}{2}}.\label{E:P}%
\end{equation}
Polynomial $P_{\omega}(t)=\big[\mathcal{R}\chi_{\Omega}\big](\omega,t)$
belongs to the range of the Radon transform. Hence it satisfies the range
conditions for Radon transform, and, in particular, the moment conditions
(see, e.g. \cite{helgason2022groups}). Namely, the $k$-th moment
\[
M_{k}(\omega)=\int\limits_{\mathbb{R}}\left[  \mathcal{R}\chi_{\Omega}\right]
(\omega,t)t^{k} \, dt=\int\limits_{h_{\Omega}^{-}(\omega)}^{h_{\Omega}^{+}%
(\omega)}P_{\omega}(t)t^{k} \, dt
\]
extends from the unit sphere $|\omega|=1$ to $\mathbb{R}^{n}$ as a homogeneous
polynomial of degree $k.$

Substituting the expression (\ref{E:P}) we have
\begin{equation}
M_{k}(\omega)=A(\omega)\int\limits_{h_{\Omega}^{-}(\omega)}^{h_{\Omega}%
^{+}(\omega)}(t-h_{\Omega}^{-}(\omega))^{\frac{n-1}{2}}(h_{\Omega}^{+}%
(\omega)-t)^{\frac{n-1}{2}}t^{k}dt. \label{E:integral}%
\end{equation}
Introduce the functions $B(\omega)$ and $C(\omega)$ as follows:
\[
C(\omega)=\frac{h_{\Omega}^{+}(\omega)-h_{\Omega}^{-}(\omega)}{2},\qquad
B(\omega)=\frac{h_{\Omega}^{+}(\omega)+h_{\Omega}^{-}(\omega)}{2}.
\]
Let us make a substitution in the integral (\ref{E:integral}):
\[
u=\frac{t-B(\omega)}{C(\omega)}.
\]
Then
\[
h_{\Omega}^{+}(\omega)=B(\omega)+C(\omega),\ h_{\Omega}^{-}(\omega
)=B(\omega)-C(\omega)
\]
and
\begin{align*}
h_{\Omega}^{+}(\omega)-t  &  =(C(\omega)+B(\omega))-(C(\omega)u+B(\omega
))=C(\omega)(1-u),\\
t-h_{\Omega}^{-}(\omega)  &  =(C(\omega)u+B(\omega))-(B(\omega)-C(\omega
))=C(\omega)(1+u).
\end{align*}
Therefore
\begin{equation}
M_{k}(\omega)=2A(\omega)C^{n}(\omega)\int_{-1}^{1}(1-u^{2})^{\frac{n-1}{2}%
}(C(\omega)u+B(\omega))^{k} \, du. \label{E:Mk}%
\end{equation}

Take $k=0:$
\[
M_{0}(\omega)=2A(\omega)C^{n}(\omega)\int_{-1}^{1}(1-u^{2})^{\frac{n-1}{2}%
} \, du.
\]
Since $M_{0}(\omega)=const,$ we obtain that the entire factor in front of the
integral is constant: $2A(\omega)C^{n}(\omega)=const=c$, so that%
\[
M_{k}(\omega)=c\int_{-1}^{1}(1-u^{2})^{\frac{n-1}{2}}(C(\omega)u+B(\omega
))^{k} \, du.
\]

Now take $k=1$:
\[
M_{1}(\omega)=c\int_{-1}^{1}(1-u^{2})^{\frac{n-1}{2}}(C(\omega)u+B(\omega
)) \, du=cB(\omega)\int_{-1}^{1}(1-u^{2})^{\frac{n-1}{2}} \, du.
\]
The first moment $M_{1}(\omega)$ extends to $\mathbb{R}^{n}$ as a linear
function, hence so does $B(\omega):$
\[
B(\omega)=b\cdot\omega+b_{0},
\]
for some vectors $b,\ b_{0}\in\mathbb{R}^{n}.$

However, $2B(\omega)=h_{\Omega}^{+}(\omega)-h_{\Omega}^{+}(-\omega)$ so that
$B(\omega)$ is an odd function of $\omega$ and, in particular, $b_{0}=0.$
Moreover, by passing to the translated domain
\[
\ \widetilde{\Omega} =\Omega-b,
\]
we can make $B(\omega)=0$ for all $\omega.$ Indeed,
\begin{align*}
h_{\widetilde{\Omega} }^{+}(\omega)  &  =\sup_{y\in\widetilde{\Omega} }%
y\cdot\omega=\sup_{x\in\Omega}(x-b)\cdot\omega=h_{\Omega}^{+}(\omega
)-b\cdot\omega,\\
h_{\widetilde{\Omega} } ^{-}(\omega)  &  =-h_{\widetilde{\Omega} }^{+}%
(-\omega)=h_{\Omega} ^{-} (\omega)-b\cdot\omega.
\end{align*}
Then for the translated domain $\widetilde{\Omega} $ we have
\[
\widetilde{B} (\omega)=\frac{1}{2}(h_{\widetilde{ \Omega}}) ^{+}(\omega) + h_{
\widetilde{ \Omega} }^{-}(\omega))=B(\omega)-b\cdot\omega=0.
\]
Thus, applying the translation $x\rightarrow x-b$ we can assume from the very
beginning that $B(\omega)=0.$ This implies $h_{\widetilde{\Omega}}%
(-\omega)=h_{\widetilde{\Omega}}(\omega)$ which means that after the
translation to the vector $b,$ domain $\Omega$ becomes centrally symmetric.
From now on, we assume that this is the case and $B=0.$

Then
\[
h_{\Omega}^{+}(\omega)=C(\omega).
\]
At last, for $k=2$ formula (\ref{E:Mk}) turns into
\[
M_{2}(\omega)=c\cdot C^{2}(\omega)\int_{-1}^{1}(1-u^{2})^{\frac{n-1}{2}}%
u^{2} \, du.
\]
Thus, $C^{2}(\omega)$ differs from $M_{2}(\omega)$ by a nonzero factor, and
since $M_{2}(\omega)$ is the restriction to $\mathbb{S}^{n-1}$ of a quadratic
homogeneous polynomial, $C^{2}(\omega)$ has the same property.

After applying an orthogonal transformation we can reduce the quadratic form
$C^{2}=(h_{\Omega}^{+})^{2}=h_{\Omega}^{2}$ to the diagonal form:
\[
h_{\Omega}^{2}(\omega)=\sum_{j=1}^{n}a_{j}\omega_{j}^{2}.
\]
Since the left hand side is strictly positive on $\mathbb{S}^{n-1}$ (indeed,
$h_{\Omega}(\omega)=0$ is impossible since it would mean $0\in\partial\Omega$
which is not the case because $\Omega$ is centrally symmetric and $0$ is its
interior point) , we have $a_{j}>0$ for all $j=1,...,n.$ We write $a_j=\alpha_j^2.$ Then we have
\[
h_{\Omega}(\omega)= |\mathcal A \omega|,
\]
where the matrix $\mathcal A $ is the non-degenerate diagonal matrix $\mathcal A =\mathrm{diag}
(\alpha_{1},...,\alpha_{n}), j=1,...,n.$

The hyperplane $x\cdot\omega=h_{\Omega}(\omega)$ is tangent to $\partial
\Omega$ and the convex domain $\Omega$ coincides with the intersection of the
open half-spaces $x\cdot\omega<h_{\Omega}(\omega)=|\mathcal A \omega|^{2},\ \omega
\in\mathbb{S}^{n-1},$ i.e.,
\[
\Omega=\{x\in\mathbb{R}^{n}:x\cdot\omega<|\mathcal A \omega|,\ \forall\omega
\in\mathbb{S}^{n-1}\}.
\]
Taking $\omega=\frac{\mathcal A ^{-1}\eta}{|\mathcal A^{-1}\eta|},$ where $\eta\in\mathbb{S}%
^{n-1}$ is arbitrary, we obtain
\[
\Omega=\{x\in\mathbb{R}^{n}:(\mathcal A ^{-1}x)\cdot\eta=x\cdot(\mathcal A^{-1}\eta)<1,\ \eta
\in\mathbb{S}^{n-1}.\}
\]
The inequality $(\mathcal A ^{-1}x)\cdot\eta<1$ for all $\eta,\ |\eta|=1,$ is equivalent
to $|\mathcal A ^{-1}x|<1$ and hence
\[
\Omega=\{x\in\mathbb{R}^{n}:|\mathcal A ^{-1}x|<1\}
\]
is bounded by the ellipsoid $\sum_{j=1}^n \frac{x_j^2}{\alpha_j^2}=1.$

\subsection{The case of even $n$}

If $k_{\Omega}=0$ and $n$ is even, then by (\ref{E:k-even})
\[
\lbrack\mathcal{H(}\mathcal{R}\chi_{\Omega})](\omega,t)=P_{\omega
}(t),\ h_{\Omega}^{-}(\omega)<t<h_{\Omega}^{+}(\omega),
\]
where $P_{\omega}$ is a polynomial of degree at most $n-1.$ By Lemma
\ref{L:Hilbert}, with $a=h_{\Omega}^{-}(\omega)$ and $b=h_{\Omega}^{+}%
(\omega)$
\[
\lbrack\mathcal{R}\chi_{\Omega}](\omega,t)=\frac{Q_{\omega}(t)}{\sqrt
{(t-h_{\Omega}^{-}(\omega))(h_{\Omega}^{+}(\omega)-t)}}.
\]
where $Q_{\omega}(t)$ is a polynomial of degree at most $n.$

By Lemma \ref{L:zeros}, we have $[\mathcal{R}\chi_{\Omega}](\omega
,t)=c(t-h_{\Omega}^{-}(\omega))^{\frac{n-1}{2}}(1+\mathit{o}(1)),$ as
$t\rightarrow h_{\Omega}^{-}(\omega)+0,$ and a similar asymptotic is true for
$t\rightarrow h_{\Omega}^{+}(\omega)-0.$ This implies that $t_{1}=h_{\Omega
}^{-}(\omega)$ and $t_{2}=h_{\Omega}^{+}(\omega)$ are zeros of the polynomial
$Q_{\omega}(t),$ each of multiplicity $\frac{n-1}{2}+\frac{1}{2}=\frac{n}{2}.$

Since, on the other hand, $\deg Q_{\omega}\leq n,$ the polynomial $Q_{\omega}$
has the representation
\[
Q_{\omega}(t)=A(\omega)(t-t_{1})^{\frac{n}{2}}(t_{2}-t)^{\frac{n}{2}}%
=A(\omega)(t-h_{\Omega}^{-}(\omega))^{\frac{n}{2}}(h_{\Omega}^{+}%
(\omega)-t)^{\frac{n}{2}}%
\]
and, correspondingly,
\[
\lbrack\mathcal{R}\chi_{\Omega}](\omega,t)=A(\omega)(t-h_{\Omega}^{-}%
(\omega))^{\frac{n-1}{2}}(h_{\Omega}^{+}(\omega)-t)^{\frac{n-1}{2}}.
\]
Then, the same argument, as in the case of odd $n,$ based on the range
description for Radon transform and using the first three moment conditions
applied to (\ref{E:integral}), implies that the support function $h_{\Omega
}=h_{\Omega}^{+}$ coincides on $\mathbb{S}^{n-1}$ with the square root of a
quadratic homogeneous polynomial, which means that $\Omega$ is an ellipsoid.
The proof of Theorem \ref{T:Main} is complete.

\section{Concluding remarks}\label{S:discussion}

Below we discuss the connections between our results and other known inversion formulas.

\begin{itemize}
\item As it was mentioned in Introduction, the problem of reconstructing a
function from its spherical means centered on a hypersurface $\Gamma$ arises
in TAT/PAT. The forward problem of TAT/PAT is modeled by the Cauchy problem
for the wave equation (see, e.g. \cite{KuKuHand})
\begin{align}
u_{tt}  &  =c^{2}(x)\Delta u,\quad t\geq0,\quad x\in\mathbb{R}^{n}%
,\label{E:wave}\\
u(x,0)  &  =f(x),\quad u_{t}(x,0)=0, \label{E:Cauchy}%
\end{align}
where $u(x,t)$ is the pressure in the propagating acoustic wave, $f(x)$ is the
initial pressure and $c(x)$\ is the speed of sound. Depending on the type of
transducers, one measures either pressure $u$ or its normal derivative
$\frac{\partial u}{\partial\nu}$ on a measurement surface $\Gamma \subset \mathbb R^n.$
The TAT/PAT inverse problem consists of  reconstructing the initial value $f(x)=u(x,0)$ from the
Dirichlet data $u|_{\Gamma\times\lbrack0,\infty)}.$  An alternative version is to find
$f(x)$ from the Neumann data $\frac{\partial}{\partial\nu}u|_{\Gamma\times
\lbrack0,\infty)}.$

For a constant speed of sound $c(x),$ solution $u(x,t)$ of
(\ref{E:wave}),\ (\ref{E:Cauchy}) can be expressed through the spherical means
$\left[  \mathcal{M}f\right]  (x,t)$ of the initial data $f(x)$ by the
Kirchoff-Poisson formula (e.g., Section 2.4, Thm 2 and 3, \cite{Evans}). This
formula allows one to reduce the problem of finding the initial value $f(x)=u(x,0)$ from the data
$u(x,t)|_{\Gamma\times \lbrack0,\infty)}$ to the problem of recovering the function $f(x)$ from
its spherical means $[\mathcal{M}f](x,t)|_{\Gamma\times \lbrack0,\infty)}$ with the centers on $\Gamma.$
The solution of the latter problem is given by formulas (\ref{E:inv-even}) and
(\ref{E:inv-odd}), and Theorem \ref{T:Main} establishes that, if
$\Gamma=\partial\Omega,$ these formulas are only valid in the case of
ellipsoidal surfaces $\partial\Omega$.

An inversion formula recovering the initial data $f(x)=u(x,0)$ from the Dirichlet data
$u|_{\Gamma\times\lbrack0,\infty)}$ was proposed in
\cite{Wang-sphere,Wang-universal,Wang-cylinder,Wang-various} for several
different acquisition surfaces $\Gamma$. Taking into account the relation
between  $u|_{\Gamma\times \lbrack0,\infty)}$  and $\mathcal{M}f|_{\Gamma\times \lbrack0,\infty)},$  one can show that the above formula is equivalent to
formulas (\ref{E:inv-even}) and (\ref{E:inv-odd}). Therefore, in the case when
$\Gamma$ is a boundary $\partial\Omega$ of a convex domain with a smooth
boundary, formula given in \cite{Wang-sphere,Wang-universal,Wang-cylinder,Wang-various} holds if and only
if $\partial\Omega$ is an ellipsoid.

\item In \cite{Haltm-Dreier-SIIMS,Haltm-Dreier-SIAM}, inversion formulas for
ellipsoids were obtained for the problem of reconstructing $f(x)$ from the
Neumann data $\frac{\partial}{\partial\nu}u|_{\Gamma\times\lbrack0,\infty)}.$
In particular, it has been shown that, as in the case of Dirichlet data, if
these formulas are applied to an arbitrary convex domain $\Omega$ with a
smooth boundary, then an error terms $K_{\Omega}$ appears, still given by
equations (\ref{E:K-kernel}), (\ref{E:k-even}), (\ref{E:k-odd}). Since we have
proven in Theorem \ref{T:Main} that $K_{\Omega}=0$ is equivalent to
$\Gamma=\partial\Omega$ being an ellipsoid, we conclude that the Neumann data
version of Theorem \ref{T:Main} is also true, i.e., inversion
formula\cite{Haltm-Dreier-SIIMS,Haltm-Dreier-SIAM} is exact for ellipsoids only.

\item Inversion formulas (\ref{E:inv-even}), (\ref{E:inv-odd}) have been also
shown \cite{Haltm-parabolic-SIAM,Haltm-universal-JMAA} to hold in the case of
some unbounded quadratic surfaces, e.g., for parabolic and elliptic cylinders.
Since an expression for the error operator is not known for the of unbounded
surfaces $\Gamma$, our analysis cannot be extended to these cases.

\item Other exact inversion formulas (formulated either in terms of the
spherical means $\mathcal{M}f|_{\Gamma\times(0,\infty)}$ or in terms of the
Dirichlet data $u|_{\Gamma\times\lbrack0,\infty)}$) have been discovered. They
hold when the surface is a sphere \cite{Finch04,Finch07} or, more generally,
an ellipsoid \cite{Salman,Pala}, and also for certain more complicated
surfaces \cite{Pala,Pal-book-2016}. These formulas are not equivalent to the
universal backprojection formulas considered here, meaning that the
corresponding inverse operators coincide only on the image of the spherical
mean operator $\mathcal{M}$ restricted to the surface $\Gamma.$ In our
opinion, it might be interesting to understand what is the largest class of
hypersurfaces $\Gamma$ for which these formulas hold.
\end{itemize}

\section{ Acknowledgement}

The question answered in this paper was posed by Professor Haltmeier in a
conversation with the first author at the 9th Conference "Inverse Problems,
Modeling and Simulation", IPMS\ 2018 held on Malta. The preparation of the
present paper started in 2022, during the 10th occurrence of this conference.
The authors thank Professor Haltmeier for the interesting question and the
organizers of IPMS-2022 for creating a stimulating environment and excellent
conditions for collaboration. The second author acknowledges support by the
NSF, through the award NSF/DMS 1814592.

\end{document}